\documentclass[12pt]{article}

\def\bbb#1{\hbox {{\gordas #1}}}
\font\gordas = msbm10 at 12pt

\def\cajita{\rule{5pt}{5pt}}

\def\bbb#1{\hbox {{\gordas #1}}}
\def\erre{{\bbb R}}

\def\ene{{\bbb N}}
\newtheorem{teorema}{Theorem}[section]
\newtheorem{proposicion}[teorema]
{Proposition}
\newtheorem{lema}[teorema]{Lemma}
\par

\par
\newtheorem{corolario}[teorema]{Corollary}
\par
\par

\par
\newtheorem{notas}[teorema]{Remarks}
\par

\par

\par

\par
\par

\begin{document}

\begin{center}
{\Large {\bf Linear Stochastic Differential Equations Driven by a Fractional
Brownian Motion with Hurst Parameter less than $1/2$}}
\end{center}

\vglue1cm
\begin{tabular}{cccc}
Jorge A. Le\'{o}n$^*$
 &  & Jaime San Mart\'{\i}n\\
Departamento de Control Autom\'atico &  & Departamento Ingenier\'{\i}a Matem\'atica\\
CINVESTAV--IPN &  & CMM,  Universidad de Chile \\
Apartado Postal 14--740 &  &Casilla 170/3\\
07000 M\'{e}xico, D.F. &  & Santiago \\
Mexico &  & Chile
\end{tabular}

\footnote{$^*$ Partially supported by CONACyT}

\begin{abstract}
In this paper we use the chaos decomposition approach to establish the existence of a unique continuous solution to linear
fractional differential equations of the Skorohod type. Here the coefficients are deterministic, the inital condition is
anticipating and the underlying fractional Brownian motion has Hurst parameter less than $1/2$. We provide an explicit
expression for the chaos decomposition of the solution in order to show our results.
\end{abstract}
\section{Introduction}
\label{sec:1}
The fractional Brownian motion (fBm) with Hurst parameter $H\in (0,1)$ is a Gaussian process with useful properties.
In particular, the stationary of its increments, and the self--similarity and the long--range dependence of this
process (see Mandelbrot and Van Ness \cite{Mv}) become the fBm a suitable driven noise for the construction of
stochastic models and the analysis of phenomena that exhibit scale--invariant and long--range correlated force.
However the fBm is not a semimartingale when $H\neq \frac{1}{2}$. Hence we cannot apply the techniques of the
stochastic calculus in the It\^o sense to define a stochastic integral with respect to the fBm.

Different interpretations of stochastic integral with respect to the fBm $B$ have been used by several authors to
study the fractional stochastic differential equation of the form
\begin{equation}
\label{eq:1.2.1}
X_t=X_0 +\int^t_0 a(s, X_s)ds+\int_0^t b(s, X_s)dB_s, \quad 0\leq t\leq T.
\end{equation}

In the case that $H\in (1/2, 1)$, it is reasonable to consider equation (\ref{eq:1.2.1}) as a path--by--path ordinary
differential equation since $B$ has H\"older--continuous paths with all exponents less than $H$ and $\int^T_0 Y_s dB_s$
exists as a pathwise Riemann--Stieltjes integral for any $\lambda$--H\"older continuous process $Y$ with $\lambda>1-H$
(see Young \cite{Yo}). This pathwise equation has been studied by several authors (see for instance \cite {Kka,
Ku, Li, Ru} and \cite{Za1}). Z\"ahle \cite{Za1} has improved this pathwise approach by constructing and extension of
the Lebesgue--Stieltjes integral via the fractional calculus \cite{Sa}. These extended integral has been considered
by Nualart and R\u{a}\c{s}canu \cite{Nr} and Z\"ahle \cite{Za2} to analyze equation (\ref{eq:1.2.1}).

The notion of $p$--variation and a limit result of Lyons \cite{Ly} allow Coutin and Qian \cite{Cq1, Cq2} to get a
Wong--Zakai type approximation limit for equation (\ref{eq:1.2.1}) when $H>1/4$.

Al\`os et al \cite{Aln} (resp. Le\'on and Tudor \cite{Lt}) work with the stochastic Stratonovich  integral in the
Russo and Vallois sense \cite{Rv} when $H\in (1/4, 1/2)$ (resp. $H\in (1/2, 1)$). Also a method based on an
extension of this Stratonovich integral is presented in \cite{Za2}.

The aim of this paper is to use the chaos expansion approach to
show an existence and uniqueness result for linear stochastic
differential equations of the form (\ref{eq:1.2.1}), in the case
that $H\in (0, 1/2)$, the coefficients are deterministic and the
stochastic integral is an extension of the divergence operator in
the Malliavin calculus sense (see Proposition \ref{P2.6.1} below).
The solution is given in terms of explicit expressions for the
kernels of the fractional multiple integrals in its chaos
expansion. Also a sufficient condition for the continuity of the
solution is provided. This chaotic expansion procedure is
introduced in Shiota \cite{Sh} for $H=1/2$. That is, when $B$ is a
Brownian motion.

The organization of the article is as follows. Section 2 describes
the framework of this paper. Namely, Section 2.1 introduces some
basic elements of the fractional calculus and Section 2.2 gives
some basic definitions and facts of the stochastic calculus for
the fBm. In Section 3 we study equation (\ref{eq:1.2.1}). Finally,
we develop the auxiliary tools needed for our proofs in Section 4.

\section{Preliminaries}
\label{sec:2}
\setcounter{equation}{0}
The purpose of this section is to describe the framework that will be used in this paper. Although some results
discussed in this section are known, we prefer to provide a self--contained exposition for the convenience of the reader.

\subsection{Fractional calculus}
\label{sub:2.1}
For a detailed account on the fractional calculus theory, we refer to Samko et al. \cite{Sa}.

Throughout $T$ is a positive number. Consider an integrable function $f:[0,T]\rightarrow \erre$ and $\alpha \in (0,1)$.
The {\it right--sided fractional integral} of $f$ of order $\alpha$ is given by
$$
I^\alpha_{T-}(f)(x)=\frac{1}{\Gamma (\alpha)}
\int^T_x \frac{f(u)}{(u-x)^{1-\alpha}}du, \quad \hbox{\rm for a.a.}\;\;x\in [0,T].
$$
Note that Fubini theorem implies that $I^\alpha_{T-}(f)$ is a
function in $L^p ([0,T]), p\geq 1$, whenever $f\in L^p ([0,T])$.
That is, the space $L^p([0,T])$ is invariant under the
right--sided fractional integral $I^\alpha_{T-}$. Actually we have
that $I^\alpha_{T-}$ has the following property (see \cite{Sa},
Theorem 3.5 and Notes \S 4):
\begin{lema}
\label{L2.1.1} {\it Let $1<p<1/\alpha$. Then $I^\alpha_{T-}$ is a linear bounded operator from $L^p ([0,T])$
into $L^r([0,T])$, for every $1\leq r\leq p(1-\alpha p)^{-1}$. }
\end{lema}

In what follows $C(p,r)$ denotes the norm of $I^\alpha_{T-}$ as a bounded linear operator from $L^p([0,T])$ into $L^r([0,T])$.

We denote by $I^\alpha_{T-}(L^p), p\geq 1$, the family of all
functions $f\in L^p ([0,T])$ such that
\begin{equation}
\label{eq:2.2.1}
f=I^\alpha_{T-}(\varphi),
\end{equation}
for some $\varphi \in L^p([0,T])$. Samko et al. \cite{Sa} (Theorem
13.2) provide a characterization of the space
$I^\alpha_{T-}(L^p)$, $p>1$. Namely, a measurable function $f$
belongs to $I^\alpha_{T-}(L^p)$ (i.e., it satisfies
(\ref{eq:2.2.1})) if and only if $f\in L^p ([0,T])$ and the
integral
\begin{equation}
\label{eq:2.2.2}
\int^T_{s+\varepsilon} \frac{f(s)-f(u)}{(u-s)^{1+\alpha}}du
\end{equation}
converges in $L^p ([0,T])$ as $\varepsilon \downarrow 0$. In this
case a function $\varphi$ satisfying (\ref{eq:2.2.1}) agrees with
the right--sided fractional derivative
\begin{equation}
\label{eq:2.2.3}
(D^\alpha_{T-}f)(s)=\frac{1}{\Gamma (1-\alpha)}
\biggl(\frac{f(s)}{(T-s)^\alpha}+\alpha \int^T_s
\frac{f(s)-f(u)}{(u-s)^{1+\alpha}}du\biggr),
\end{equation}
where the integral is the $L^p([0,T])$--limit of (\ref{eq:2.2.2}). Moreover, by \cite{Sa} (Lemma 2.5),
there is at most one solution $\varphi$ to the equation (\ref{eq:2.2.1}).

A useful tool to analyze the convergence of (\ref{eq:2.2.2}) is
the following inequality.
\begin{lema}
\label{L2.2.1} Let $0<s<t<r$. Then
$$
\alpha \int^r_t (r-u)^{\alpha -1} u^{-\alpha}(u-s)^{-\alpha -1}du
\leq t^{-\alpha}(t-s)^{-\alpha} (r-s)^{-1} (r-t)^\alpha.
$$
\end{lema}
{\it Proof:} The result is an immediate consequence of the changes
of variables $z=(r-u)(u-s)^{-1}$ and $y=(t-s)z$, and the fact that
$y\mapsto \frac{(t-s)r+ys}{y+t-s}$ is a decreasing function on
$[0, r-t]$. \hfill $\cajita$

In the remaining of this paper, we will consider the space
\begin{equation}
\label{eq:2.2.4} \Lambda^\alpha_T:=\{f:[0,T]\rightarrow \erre :
\exists \;\phi_f \in L^2 ([0,T]), \hbox{ such that }  f(u)
=u^\alpha I^\alpha_{T-}(s^{-\alpha} \phi_f (s))(u)\}.
\end{equation}
Here $\alpha \in (0, \frac{1}{2})$. More precisely, we will make use of the following result.
\begin{lema}
\label{L2.3.1} Let $f\in \Lambda^\alpha_T$ be such that $\phi_f \in L^p ([0,T])$ for some $p\in (2, 1/\alpha)$. Then $f1_{[0,t]}$
also belongs to $\Lambda^\alpha_T$ for each $t\in [0,T]$, and for any $p'\in [2,p)$,
$$
||\phi_{f1_{[0,t]}}||_{L^{p'}([0,T])}\leq C_{\alpha, p', p, t}||\phi_f||_{L^p([0,T])},
$$
where $C_{\alpha, p', p, t}=t^{(p-p')/p'p} + \frac{C(p, p(1-\alpha
p)^{-1})}{\Gamma (1-\alpha)}(\int^t_0 (t-s)^{-p'\alpha
q}ds)^{1/p'q} $ with $q=\frac{p}{p-p'(1-\alpha p)}.$
\end{lema}
\noindent
{\bf Remarks.}
\begin{itemize}
\item[i)] The assumptions of the result give that $1-p'\alpha q>0$.
\item[ii)] Remember that $C(p, p(1-\alpha p)^{-1})$ is the norm of the linear operator
$I^\alpha_{T-}:L^p ([0,T])\rightarrow L^{p(1-\alpha p)^{-1}}([0,T])$ (see Lemma \ref{L2.1.1}).\end{itemize}
{\it Proof:} Fix $t\in [0,T]$. Then, by (\ref{eq:2.2.3}), we have
\begin{eqnarray*}
(D^\alpha_{T-} (u^{-\alpha}f(u)1_{[0,t]}(u)))(s)
&=& 1_{[0,t]}(s) \biggl[s^{-\alpha}\phi_f (s)\\
&&+ \frac{\alpha}{\Gamma (1-\alpha)}\int^T_t
\frac{u^{-\alpha }f(u)}{(u-s)^{1+\alpha}}du\biggr].
\end{eqnarray*}
That is (see (\ref{eq:2.2.4})),
\begin{equation}
\label{eq:2.3.1}
\phi_{f1_{[0,t]}}(s)=1_{[0,t]}(s)
\biggl[\phi_f(s) +\frac{\alpha s^\alpha}{\Gamma (1-\alpha)}
\int^T_t\frac{u^{-\alpha}f(u)}{(u-s)^{1+\alpha}}du\biggr].
\end{equation}

Finally, observe that Lemma \ref{L2.2.1} gives
\begin{eqnarray*}
\lefteqn{
||(\cdot)^\alpha 1_{[0,t]}(\cdot)\int^T_t
\frac{u^{-\alpha}|f(u)|}{(u-\cdot)^{1+\alpha}}
du||_{L^{p'}([0,T])}}\\
&\leq & \frac{1}{\alpha}||1_{[0,t]} (\cdot ) (t-\cdot)^{-\alpha}
\frac{1}{\Gamma (\alpha)}\int^T_t \frac{|\phi_f
(r)|}{(r-\cdot)^{1-\alpha}} dr||_{L^{p'}([0,T])}.
\end{eqnarray*}
Hence, Lemma \ref{L2.1.1} implies that for $q=\frac{p}{p-p'(1-\alpha p)}$,
\begin{eqnarray*}
\lefteqn{
||(\cdot)^\alpha 1_{[0,t]}(\cdot )\int^T_t
\frac{u^{-\alpha}|f(u)|}{(u-\cdot)^{1+\alpha}}
du||_{L^{p'}([0,T])}}\\
&\leq & \frac{C(p, p(1-\alpha p)^{-1})}{\alpha} (\int^t_0
(t-s)^{-p'\alpha q}ds)^{1/p'q}||\phi_f||_{L^p ([0,T])},
\end{eqnarray*}
which, together with (\ref{eq:2.3.1}), yields that the Lemma holds. \hfill $\cajita$

We will also need the following result.
\begin{lema}
\label{L2.4.1} Let $f$ be a function in $\Lambda^\alpha_T$ and $g:[0,T]\rightarrow \erre$ a H\"older continuous function
with parameter $\beta >\alpha$. Then $gf$ also belongs to $\Lambda^\alpha_T$.
\end{lema}
\noindent
{\it Proof:} Using equality (\ref{eq:2.2.3}) again, we obtain
\begin{eqnarray*}
(D^\alpha_{T-}(u^{-\alpha}g(u) f(u)))(s) &=& g(s) (D^\alpha_{T-}(u^{-\alpha}f(u)))(s)\\
&&+ \frac{\alpha}{\Gamma (1-\alpha)}\int^T_su^{-\alpha} f(u)
\frac{g(s)-g(u)}{(u-s)^{\alpha +1}}du.
\end{eqnarray*}
Finally, it follows from the H\"older continuity of $g$ that $s\mapsto s^\alpha
\int^T_s u^{-\alpha} |f(u)|\frac{|g(s)-g(u)|}{(u-s)^{\alpha +1}}du$ is a square--integrable function. Thus $gf$ belongs to $\Lambda^\alpha_T$. \hfill $\cajita$
\subsection{Fractional Brownian motion}
\label{sub:2.3}
Throughout $B^H=\{B^H_t:t\in [0,T]\}$ is a fractional Brownian motion (fBm) with Hurst parameter $H\in (0, 1/2)$
defined on a complete probability space $(\Omega, {\cal F}, P)$. It means, the fBm $B^H$ is a Gaussian process with zero mean and covariance function
$$
R_H(t,s):=\frac{1}{2}(t^{2H}+ s^{2H}-|t-s|^{2H}), \quad s, t\in [0,T].
$$
In the remaining of this paper, we assume ${\cal F}=\sigma \{B_t:t\in [0,T]\}$.
The reader can consult Nualart \cite{Dn} and references therein for a recent presentation of the facts related to the fBm.

Let ${\cal H}_H$ be the Hilbert space defined as the completion of  the step functions on $[0,T]$ with respect to the inner product
$$
\langle 1_{[0,t]}, 1_{[0,s]}\rangle_{{\cal H}_H}=R_H(t,s)\quad t, s\in [0,T].
$$
From Pipiras and Taqqu \cite{Pt} (see also \cite{Dn}), it follows that ${\cal H}_H$ coincides with the Hilbert space
$\Lambda^{1/2-H}_T$ (introduced in (\ref{eq:2.2.4})) equipped with scalar product
$$
\langle f,g\rangle_{\Lambda_T^{1/2 -H}}=C_H \langle \phi_f, \phi_g\rangle_{L^2([0,T])},
$$
with $C_H =\frac{2H\Gamma (H+\frac{1}{2})}{(1-2H)\beta (1-2H, H +\frac{1}{2})}$. So the map $1_{[0,t]}\mapsto B^H_t$ is
extended to an isometry of $\Lambda^{1/2-H}_{T}$ onto a Guassian closed subspace of $L^2(\Omega, {\cal F}, P)$
(see Nualart \cite{Dn}). This isometry is denoted by $\phi\mapsto B^H(\phi)$.

Now we assume that the reader is familiar with the basic elements of the stochastic calculus of variations for
Gaussian processes as given for example in Nualart \cite{Nu}.

Let $n$ be a positive integer. The $n$-th multiple integral $I_n$ of order $n$ with respect to $B^H$ is a linear
operator from the $n$--th symmetric tensor product $(\Lambda^{1/2 -H}_T)^{\odot n}$ of $\Lambda^{1/2-H}_T$ into
$L^2(\Omega, {\cal F}, P)$ satisfying the following two properties:
\begin{itemize}
\item[$\bullet$] Let $H_m$ be the $m$--th Hermite polynomial
$$
H_m(x)=\frac{(-1)^m}{m!}e^{x^2/2} \frac{d^m}{dx^m} e^{-x^2/2}, \quad x\in \erre,$$
and $\{e_k:k\in \ene\}$ an orthonormal system on $\Lambda^{1/2-H}_T$. Then, for any $f_n \in (\Lambda^{1/2-H})^{\odot n}$, we have
\begin{eqnarray*}
\lefteqn{
E[I_n (f_n)(n_1 !)H_{n_1}(B^H(e_{i_1}))\cdots (n_k )!H_{n_k}(B^H(e_{i_k}))]}\\
&=& \left\{
\begin{array}{lcl}
n!\langle f_n, e^{\otimes n_1}_{i_1}\otimes \cdots \otimes e^{\otimes n_k}_{i_k}\rangle _{(\Lambda^{1/2 -H}_T)^{\otimes n}}, & \hbox{\rm if} & n=\sum^k_{j=1}n_j\\
0, && \hbox{\rm otherwise}.
\end{array}\right.
\end{eqnarray*}
\item[$\bullet$] Let $f\in (\Lambda^{1/2-H}_T)^{\odot m}$ and $g\in (\Lambda^{1/2-H}_T)^{\odot n}$. Then
$$E[I_m (f) I_n (g) ]=\left\{
\begin{array}{lcl}
0, & \hbox{\rm if} & n\neq m,\\
m!\langle f, g\rangle_{(\Lambda^{1/2 -H}_{T})^{\otimes m}}, & \hbox{\rm if} & n=m.
\end{array}\right.
$$
\end{itemize}
As a consequence of the relation between multiple integrals and Hermite polynomials, we have that
any $F\in L^2 (\Omega, {\cal F}, P)$ has a unique chaotic representation of the form
$$
F=\sum^\infty_{n=0} I_n (f_n),
$$
with $I_0 (f_0)=EF$.

Le\'on and Nualart \cite{Ln} have extended the domain of the divergence operator in the sense of Malliavin calculus
for Gaussian stochastic processes. This extension was first analyzed by Cheridito and Nualart \cite{Cn}
when the underlying Gaussian process is the fBm $B=\{B_t:t\in \erre\}$ with Hurst parameter $H\in (0, 1/2)$.
For the fBm $B^H$, this extension denoted by $\delta$ is characterized by the following result (see \cite{Ln}).
\begin{proposicion}
\label{P2.6.1} Let $u\in L^2 (\Omega; L^2 ([0,T]))$ be a random
variable with the chaos representation
$$
u=\sum^\infty_{n=0}I_n (f_n),\;\;\; f_n \in (\Lambda^{1/2-H}_T)^{\odot n}\otimes L^2 ([0,T]).
$$
Then $u\in\;\hbox{\rm Dom}\;\delta$ iff $\widetilde{f}_n$ (the symmetrization of $f_n$ in $L^2 ([0,T]^{n+1}))$ is in $(\Lambda^{1/2-H}_T)^{\odot (n+1)}$ for every $n\in \ene$ and
$$
\sum^\infty_{n=1} n! ||\widetilde{f}_{n-1}||_{(\Lambda^{1/2-H}_T)^{\otimes n}} <\infty .
$$
In this case $\delta (u) =\sum\limits^\infty_{n=1}I_n (\widetilde{f}_{n-1})$.
\end{proposicion}
\noindent
{\bf Remarks.}
\begin{itemize}
\item[i)] The space $\Lambda^\alpha_T$ is included in $L^2([0,T])$ for any $\alpha\in (0, 1/2)$, by the Fubini's theorem.
\item[ii)] In \cite{Cn} and \cite{Ln}, it is shown that the domain of $\delta$ is bigger than that of
the usual divergence operator. Also Hu \cite{Hu} (Section 7.1) has considered a set of integrable
processes including Dom $\delta$. The reader can see Decreusefond and \"Ust\"unel \cite{Ud} for a related
construction of a stochastic integral with respect to $B^H$.
\item[iii)] In Section 3, we use the convention
$$
\int^t_0 u_s dB_s^H=\delta (u 1_{[0,t]}),
$$
whenever $u1_{[0,t]}\in$ Dom $\delta$.
\end{itemize}
\section{Linear fractional differential equations}
\label{sec:3}
\setcounter{equation}{0}

In the remaining of this paper $B^H$ is a fBm with Hurst parameter $H\in (0,1/2)$ and we use the notation
$$
\alpha =\frac{1}{2}-H.
$$

In this section we study the chaos decomposition of the solution to a linear stochastic differential equation of the form
\begin{equation}
\label{eq:3.1.1}
X_t =\eta +\int^t_0 a(s) X_s ds+\int^t_0 b(s)X_s dB^H_s, \quad t\in [0,T].
\end{equation}
Here $\eta$ is a square--integrable random variable having the
chaotic representation
\begin{equation}
\label{eq:3.1.2}
\eta =\sum^\infty_{n=0}I_n (\eta_n),
\end{equation}
$a,b$ are two functions in $L^2([0,T])$ such that $b$ is also in $\Lambda^\alpha_T$.
\subsection{Statement of problem and main results}
\label{sub:3.1}
Suppose that equation (\ref{eq:3.1.1}) has a solution $X$ in $L^2(\Omega \times [0,T])$ with the chaos decomposition
\begin{equation}
\label{eq:3.1.3}
X_t =\sum^\infty_{n=0} I_n (f^t_n),\quad f_n \in (\Lambda^\alpha_T)^{\odot n}\otimes L^2 ([0,T]).
\end{equation}
Then the uniqueness of the chaotic representation (\ref{eq:3.1.3}), and Proposition \ref{P2.6.1} imply
$$
f^t_0 =\eta_0 \exp (\int^t_0 a(s) ds),\quad t\in [0,T],
$$
and
\begin{eqnarray*}
f^t_n(t_1, \ldots, t_n) &=& \eta_n (t_1, \ldots , t_n) +\int^t_0 a(s) f^s_n (t_1, \ldots, t_n) ds\\
&&\quad +\frac{1}{n}\sum^n_{j=1}b(t_j) f^{t_j}_{n-1} (t_1, \ldots ,\hat{t}_j, \ldots, t_n)1_{[0,t]}(t_j),\;\;\; t, t_1, \ldots, t_n \in [0,T].
\end{eqnarray*}
Hence, using induction on $n$, we obtain
\begin{eqnarray}
\label{eq:3.2.1}
\lefteqn{f_n^t (t_1, \ldots, t_n) =\exp(\int^t_0 a(s) ds)\biggl[\eta_n (t_1, \ldots, t_n)}\nonumber\\
&+& \sum^n_{j=1}\sum_{\Delta_{j,n}}\frac{(n-j)!}{j!n!}b^{\otimes j}(t_{i_1},\ldots, t_{i_j})\eta_{n-j}  (\hat{t}_{i_{1}}, \ldots, \hat{t}_{i_j})1_{[0,t]^j}
(t_{i_1},\ldots , t_{i_j})\biggr].
\end{eqnarray}
Here we use the convention
$$
\Delta_{j,n}=\{\{i_1, \ldots,i_j\}\subset \{1,\ldots, n\}:i_k \neq i_\ell \;\;\hbox{\rm if} \;\;k\neq \ell\}.
$$
Consequently, equation (\ref{eq:3.1.1}) has at most one solution in $L^2(\Omega \times [0,T])$.

Conversely, suppose that the functions $f_n$ given by (\ref{eq:3.2.1}) satisfy the following conditions:
\begin{enumerate}
\item For every $n\geq 0$, $f_n \in (\Lambda^\alpha_T)^{\odot n}\otimes L^2 ([0,T])$.
\item The process $Y_t=\sum\limits^\infty_{n=0} I_n (f^t_n)$ belongs to $L^2(\Omega \times [0,T])$. That is,
\begin{equation}
\label{eq:3.2.2}
\sum^\infty_{n=0} n! \int^T_0 ||f^t_n ||^2_{(\Lambda^\alpha_T)^{\otimes n}}dt<\infty.
\end{equation}
\item For almost all  $t\in [0,T]$, $bY 1_{[0,t]}$ belongs to Dom $\delta$.
\end{enumerate}
Then the process $Y$ is a solution in $L^2(\Omega \times [0,T])$ of equation (\ref{eq:3.1.1}).

The approach just described allows us to show the following result, which is one our goals in this paper.

Henceforth we use the notation
$$
B_{H,p,t}=1+C_H \biggl(B_{T,p}||\phi_{b1_{[0,t]}}||_{L^p([0,T])}+||\phi_{b1_{[0,t]}}||_{L^2([0,T])}\biggr)^2
$$
with
$$
B_{T,p}=\frac{T^{(p-2)/2p}}{\Gamma (1-\alpha)}C(p, p/(1-\alpha
p))\left(\frac{p-2(1-\alpha p)}{p-2}\right)^{(p-2(1-p\alpha))/2p}.
$$
\begin{teorema}
\label{T3.3.1} Let $p\in (2, 1/\alpha)$. Assume that $\phi_b \in L^p([0,T])$ and
\begin{equation}
\label{eq:3.3.1}
\sum^\infty_{k=0} (k+1)!||\eta_k||^2_{(\Lambda ^\alpha_T)^{\otimes k}}
(\sup\limits_{t\in [0,T)}B_{H,\tilde{p},t})^k <\infty
\end{equation}
for some $\tilde{p}\in (2,p)$. Then equation (\ref{eq:3.1.1}) has a unique solution in $L^2 (\Omega \times [0,T])$ given by
$$
X_t=\sum^\infty_{n=0} I_n (f^t_n),
$$
where $f_n$ is defined in (\ref{eq:3.2.1}).
\end{teorema}

\begin{notas}
\label{N3.3.2}

\begin{itemize}
\item[i)] {\rm In Lemma \ref{L2.3.1} we have found a  bound for
$||\phi_{b1_{[0,t]}}||_{L^{\tilde{p}}([0,T])}$, ${\tilde{p}}\in [2,p)$.}
\item[ii)] {\rm The following are examples of initial conditions that satisfy Hypothesis (\ref{eq:3.3.1}):}
\begin{itemize}
\item[a)] {\rm $\eta$ has a finite chaos decomposition. That is
$\eta=\sum\limits^M_{n=0} I_n (\eta_n),\  M<\infty $.}
\item[b)] {\rm $\eta$ has exponential growth. It means, there is a positive constant $C$ such that
$$
||\eta_n||_{(\Lambda^\alpha_T)^{\otimes n\leq}} C^n(n!)^{-1}.
$$}
\item[c)] {\rm There is $\varepsilon >0$ such that
$$
\sum^\infty_{k=0}(k!)^{1+\varepsilon}||\eta_k||^2_{(\Lambda^\alpha_T)^{\otimes n}}<\infty.
$$}
\end{itemize}
\item[iii)] {\rm By (\ref{eq:3.2.1}) we have that the solution $X$ of equation (\ref{eq:3.1.1}) has the form
$$
X_t =\exp (\int^t_0 a(s)ds)Y_t,
$$
where $Y$ is the solution in $L^2(\Omega \times [0,T])$ to the equation
$$
Y_t =\eta +\int^t_0 b(s)Y_s dB_s, \quad t\in [0,T].
$$}
\end{itemize}
\end{notas}

The Kolmogorov's continuity criterion implies the following result, which is our second goal in this paper.

\begin{teorema}
\label{T3.4.1} Let $p, \phi_b$ and $\eta$ be as in Theorem \ref{T3.3.1}. Moreover assume
\begin{equation}
\label{eq:3.4.1}
\sum^\infty_{k=1} e^{k\theta}k^{k/2}||\eta_k ||_{(\Lambda ^\alpha_T)^{\otimes k}}<\infty,
\end{equation}
for some $\theta $ such that $(1+e^{2\theta})(\frac{p-2}{2p}\wedge \alpha \wedge \frac{(1-\alpha p)}{p})>1$.
Then the solution of equation (\ref{eq:3.1.1}) has a continuous version in $L^2 (\Omega \times [0,T])$.
\end{teorema}
\noindent
{\bf Remarks.}
\begin{itemize}
\item[i)]Observe that (\ref{eq:3.3.1}) holds (resp. (\ref{eq:3.4.1}) is not
 true) when
 $$||\eta_k ||_{(\Lambda^\alpha_T)^{\otimes k}}=(k^2(k+1)!
 (\sup\limits_{t\in [0,T)}B_{H, \tilde{p},t})^k)^{-1/2}$$
 (resp. and $e^\theta \geq (\sup\limits_{t\in [0,T)}B_{H, \tilde{p},t}))$.
\item[ii)] Condition (\ref{eq:3.4.1}) does not necessarily imply
Assumption (\ref{eq:3.3.1}) when $\theta$ is
such that $e^{\theta +\varepsilon}\leq (\sup\limits_{t\in [0,T)}B_{H, \tilde{p},t})^{1/2}$
for some $\varepsilon >1/2$. Indeed, if we have
$$
||\eta_{k}||_{(\Lambda^\alpha_T)^{\otimes k}}= (\sup\limits_{t\in
[0,T)}B_{H, \tilde{p},t})^{-k/2}((k+1)!)^{-1/2}.
$$
Then(\ref{eq:3.4.1}) is satisfied but (\ref{eq:3.3.1}) does not
hold.
\item[iii)] Remark \ref{N3.3.2} ii) is an example of initial conditions
that satisfy (\ref{eq:3.3.1}) and (\ref{eq:3.4.1}) at the same time.
\end{itemize}

\subsection{Proof of the main results}
\label{sub:3.2}
We begin this section with the proof of Theorem \ref{T3.3.1}.
\vglue .25cm
\noindent
{\it Proof of Theorem \ref{T3.3.1}:} We will follow the method (Steps 1--3) indicated in Section 3.1.

Let $f_n$ be given by (\ref{eq:3.2.1}). Then for each $t\in [0,T]$,
\begin{eqnarray}
\label{eq:3.6.1}
||f^t_n ||_{(\Lambda^\alpha_T)^{\otimes n}} &\leq & \exp (\int^t_0 a(s) ds)\biggl[||\eta_n ||_{(\Lambda^\alpha_T)^{\otimes n}}\nonumber\\
&&\quad + \sum^n_{j=1}(j!)^{-1}||b1_{[0,t]}||^j_{\Lambda^\alpha_T}||\eta_{n-j}||_{(\Lambda^\alpha_T)^{\otimes (n-j)}}\biggr]\nonumber\\
&=& \exp (\int^t_0 a(s)ds)\sum^n_{j=0}\frac{||b1_{[0,t]}||^{n-j}_{\Lambda^\alpha_T}}{(n-j)!}||\eta_j||_{(\Lambda^\alpha_T)^{\otimes j}}.
\end{eqnarray}
Thus Lemma \ref{L2.3.1}
yields that $f_n \in (\Lambda^\alpha_T)^{\odot n}\otimes L^2 ([0,T])$.

Now we see that (\ref{eq:3.2.2}) is satisfied. By (\ref{eq:3.6.1}) we have
\begin{eqnarray*}
\lefteqn{\sum^\infty_{n=0} n!||f^t_n ||^2_{(\Lambda^\alpha_T)^{\otimes n}}}\\
&\leq & \exp (2\int^t_0 a(s) ds) \sum^\infty_{n=0} n!\left(\sum^n_{j=0}
\frac{||b1_{[0,t]}||^{2(n-j)}_{\Lambda^\alpha_T}}{j!(n-j)!}\right)\left(\sum^n_{j=0}j!\frac{||\eta_j ||^2_{(\Lambda^\alpha _T)^{\otimes j}}}{(n-j)!}\right)\\
&=& \exp (2 \int^t_0 a(s) ds)\sum^\infty_{n=0} (||b1_{[0,t]}||^2_{\Lambda^\alpha_T}+1)^n \sum^n_{j=0}j!
\frac{||\eta_j||^2_{(\Lambda^\alpha_T)^{\otimes j}}}{(n-j)!}\\
&=&\exp (2 \int^t_0 a(s) ds+||b1_{[0,t]}||^2_{\Lambda^\alpha_T}+1)\sum^\infty_{j=0} j!||\eta_j ||^2_{(\Lambda ^\alpha_T)^{\otimes j}}
(||b1_{[0,t]}||^2_{\Lambda^\alpha_T}+1)^j.
\end{eqnarray*}
Consequently, from (\ref{eq:3.3.1}) and Lemma \ref{L2.3.1} it follows that (\ref{eq:3.2.2}) holds and
\begin{equation}
\label{eq:n50}
\sum^\infty_{n=0}n!\int^T_0b(t)^2 ||f^t_n||^2_{(\Lambda^\alpha_T)^{\otimes n}}dt<\infty.
\end{equation}
Set $Y_t =\sum\limits^\infty_{n=0} I_n (f^t_n)$. So, to finish the proof, we only need to show that for each $t\in [0,T]$, the process $bY 1_{[0,t]}$ belongs to Dom $\delta$.

By Lemmas \ref{L4.5.2}, \ref{L4.7.1} and \ref{L4.8.1} below we have
\begin{eqnarray*}
&&\sum^\infty_{n=0}
(n+1)!||\frac{1}{n+1}\sum^{n+1}_{k=1}1_{[0,t]}(t_k)
b(t_k)f^{t_k}_n(\hat{t}_k)||^2_{(\Lambda^\alpha_T)^{\otimes (n+1)}}\leq\\
& &2 \exp\left(2\int_0^T |a(s)| ds + 1 + C_H
||\phi_{b1_{[0,t]}}||^2_{L^2([0,T])}\right)\sum^\infty_{k=0}k!||\eta_k||^2_{(\Lambda^\alpha_T)^{\otimes
k}}(1+C_H||\phi_{b1_{[0,t]}}
||^2_{L^2([0,T])})^k\\
& &+2\biggl(\frac{A\alpha}{\Gamma (1-\alpha)}\biggr)^2C_H \exp
(B_{H, \tilde{p}, t}) \biggl\{\sum^\infty_{k=0}(k+1)!||\eta_k
||^2_{(\Lambda^\alpha_T)^{\otimes k}} (B_{H, \tilde{p}, t})^k
+\sum^\infty_{k=0}k!||\eta_k||^2_{(\Lambda^\alpha_t)^{\otimes
k}}(B_{H, \tilde{p}, t})^{k+1}\biggr\}.
\end{eqnarray*}
Thus the result follows from (\ref{eq:3.3.1}), (\ref{eq:n50}) and Proposition \ref{P2.6.1}. \hfill $\cajita$

\medskip

Now we give the proof of Theorem \ref{T3.4.1}.

\noindent
{\it Proof of Theorem \ref{T3.4.1}:} We first note that we can assume that $a=0$ using Remark \ref{N3.3.2}. iii).
So in this case, by the Hypercontractivity property (see \cite{Nu}) and (\ref{eq:3.2.1}), the solution $X$ of equation (\ref{eq:3.1.1}) satisfies for $q=1+e^{2\theta}$,
\begin{eqnarray}
\label{eq:3.8.1}
E|X_t-X_s|^q &\leq & [\sum^\infty_{n=1}||I_n (f_n^t-f^s_n)||_{L^q(\Omega)}]^q\nonumber\\
&\leq & (t-s)^{\delta q}[\sum^\infty_{n=1}e^{n\theta}\sqrt{n!}\sum^n_{j=1}\frac{1}{j!}||\eta_{n-j}||_{(\Lambda^\alpha_T)^{\otimes (n-j)}}(C||\phi_b||_{L^p([0,T])})^j]^q,
\end{eqnarray}
where the last inequality follows from Corollary \ref{C4.12.1} below.

Now observe that there exists a positive constant $\tilde{C}$ such that
\begin{eqnarray}
\label{eq:3.9.1}
\lefteqn{
\sum^\infty_{n=1}e^{n\theta}\sqrt{n!}\sum^n_{j=1}\frac{1}{j!}||\eta_{n-j}
||_{(\Lambda^\alpha_T)^{\otimes (n-j)}}(C||\phi_b||_{L^p([0,T])})^j}\nonumber\\
&=&\sum^\infty_{k=0}||\eta_k||_{(\Lambda^\alpha_T)^{\otimes k}}\sum^\infty_{n=k+1}
\frac{\sqrt{n!}e^{n\theta}}{(n-k)!}(C||\phi_b||_{L^p([0,T])})^{n-k}\nonumber\\
&\leq &\sum^\infty_{k=0}e^{k\theta}||\eta_k||_{(\Lambda^\alpha_T)^{\otimes k}}\sum^\infty_{n=1}\sqrt{(n+k)^k}
\frac{(e^\theta C||\phi_b||_{L^p ([0,T])})^n}
{\sqrt{n!}}\nonumber\\
&\leq & \sum^\infty_{k=0}e^{k\theta} ||\eta_k||_{(\Lambda^\alpha_T)^{\otimes k}}k^{k/2} \sum^\infty_{n=1}
\frac{(e^{\theta +\frac{1}{2}}C||\phi_b||_{L^p ([0,T])})^n}
{\sqrt{n!}}\nonumber\\
&\leq& \tilde{C}\sum^\infty_{k=0} e^{k\theta}k^{k/2} ||\eta_k||_{(\Lambda^\alpha_T)^{\otimes k}}.
\end{eqnarray}

Finally, the Kolmogorov's continuity theorem, together with (\ref{eq:3.8.1}) and (\ref{eq:3.9.1}), completes the proof. \hfill $\cajita$

\section{Appendix}
\setcounter{equation}{0}
In this section the basic tool for our proofs is established.

Let $x\in \erre^n$ and $(r_{i_1}, \ldots r_{i_j})\in \erre^j$, $j\leq n$. Henceforth we use the notation
\begin{eqnarray*}
\Delta^x_k (r_{i_1},\ldots , r_{i_j})&=&\bigg\{u\in \erre^n:u_i
=x_i\;\;\hbox{\rm for}
\;\;i\not\in \{i_1, \ldots, i_j\},\\
&&u_i\in\{r_i, x_i\}\;\;\hbox{\rm for}\;\;i\in \{i_1, \ldots, i_j\} \;\;\hbox{\rm and}\\
&&\#\{u_{i_1},\ldots , u_{i_j}\}\cap \{r_{i_1},\ldots,
r_{i_j}\}=k\bigg\},
\end{eqnarray*}
$S_{j,n}$ for the class of ordered subsets
$$
\{i_1<\cdots <i_j\}\subset \Delta_{j,n}=\{\{i_1,
\ldots,i_j\}\subset \{1,\ldots, n\}:i_k \neq i_\ell \;\;\hbox{\rm
if} \;\;k\neq \ell\},
$$
and $e(t)=\exp (\int^t_0 a(s)ds)$.

\begin{lema}
\label{L4.1.1} Let $f\in \Lambda^\alpha_T$. Then for every $n\in \ene$,
$(t_1, \ldots , t_n)\mapsto f^{\otimes n}(t_1, \ldots, t_n)e (t_1 \vee \cdots \vee t_n)\in (\Lambda^\alpha_T)^{\otimes n}$ and
\begin{eqnarray*}
\lefteqn{
f^{\otimes n}(t_1, \ldots, t_n)
e(t_1\vee \cdots \vee t_n) (t_1 \cdots t_n)^{-\alpha}}\\
&=&\Gamma (\alpha)^{-n}\biggl\{\int^T_{t_1}\cdots \int^T_{t_n}
\frac{
(\prod\limits^{n}_{i=1}\phi_f(s_i)s_i^{-\alpha})e(s_1\vee \cdots \vee s_n)}
{\prod\limits^n_{i=1}(s_i-t_i)^{1-\alpha}}
ds_n\cdots ds_1\\
&&+\sum^n_{j=1}(\Gamma (1-\alpha))^{-j}\alpha^j \sum_{S_{j,n}}\int^T_{t_1}\cdots \int^T_{t_n}
\frac{
\prod_{i\not\in \{i_1, \ldots, i_j\}}\phi_f (s_i)(s_i)^{-\alpha}}
{\prod\limits^n_{i=1}(s_i-t_i)^{1-\alpha}}\\
&&\quad \cdot \biggl[\int^T_{s_{i_1}}\cdots \int^T_{s_{i_j}}
\frac{ \prod\limits^j_{k=1}f(u_{i_k})(u_{i_k})^{-\alpha}}
{\prod\limits^j_{k=1}(u_{i_k}-s_{i_k})^{1+\alpha}}
\sum^j_{k=0}(-1)^k \!\!\!\!\!\!\sum_{\delta \in \Delta^s_k
(u_{i_1},\ldots, u_{i_j})} \!\!\!\!\!\!e(\delta_1 \vee \cdots \vee
\delta_n)du_{i_j}\cdots du_{i_1}\biggr] ds_n \cdots ds_1\biggr\}.
\end{eqnarray*}
\end{lema}
\noindent
{\bf Remark.} In Lemma \ref{L4.1.1} we are using the convention
 $\prod_{i\not\in S_{n,n}}\phi_f (s_i)(s_i)^{-\alpha } =1$.
\vglue .5cm
\noindent
{\it Proof:} The result holds for $n=1$ due to Lemma \ref{L2.4.1}.
Now we use induction on $n$. So we assume that the result is true for $n$.
Then, Lemma \ref{L2.4.1} implies
\begin{eqnarray*}
&&f(t_{n+1})t^{-\alpha}_{n+1}((\prod\limits^{n}_{i=1}f(t_i)t^{-\alpha}_i)e (t_1\vee \cdots \vee t_{n+1}))\\
&=& (\Gamma (\alpha))^{-(n+1)}\biggl\{\int^T_{t_1}\cdots \int^T_{t_{n+1}}\frac{
\prod\limits^n_{i=1}\phi_f (s_i)s^{-\alpha}_i}
{\prod\limits^{n+1}_{i=1}(s_i-t_i)^{1-\alpha}}
\biggl[e(s_1\vee \cdots \vee s_{n+1})\phi_f (s_{n+1})s^{-\alpha}_{n+1}\\
&& + \frac{\alpha}{\Gamma(1-\alpha)} \int^T_{s_{n+1}}
\phi_f(u)u^{-\alpha} \frac{e(s_1\vee \cdots \vee
s_{n+1})-e(s_1\vee \cdots \vee s_n \vee u)}
{(u-s_{n+1})^{1+\alpha}}du
\biggr] ds_{n+1}\cdots ds_1\biggr\}\\
&& +(\Gamma
(\alpha))^{-n}\sum^n_{j=1}\left(\frac{\alpha}{\Gamma(1-\alpha)}\right)^j
\sum_{S_{j,n}}\int^T_{t_1}\cdots \int^T_{t_n}
\frac{\prod\limits_{i\not\in \{i_1,\ldots , i_j, n+1\}}\phi_f
(s_i)s_i^{-\alpha}}
{\prod\limits^n_{i=1}(s_i-t_i)^{1-\alpha}}\\
&&\cdot \int^T_{s_{i_1}}\cdots  \int^T_{s_{i_j}} \frac{
\prod\limits^j_{k=1}f(u_{i_k})u^{-\alpha}_{i_k}}
{\prod\limits^j_{k=1}(u_{i_k}-s_{i_k})^{1+\alpha}}
\sum^j_{k=0}(-1)^k \sum_{\delta \in \Delta^{(s_1, \ldots ,s_n)}_k (u_{i_1}, \ldots , u_{i_j})}\\
&&\cdot \left[\frac{1}{\Gamma (\alpha)}\int^T_{t_{n+1}} \frac{
\phi_f (s_{n+1})s^{-\alpha}_{n+1}e(\delta_1 \vee\cdots \vee
\delta_n \vee s_{n+1})} {(s_{n+1}-t_{n+1})^{1-\alpha}}
ds_{n+1}\right.\\
&& + \frac{\alpha}{\Gamma (\alpha) \Gamma (1-\alpha)}
\int^T_{t_{n+1}}(s_{n+1}-t_{n+1})^{\alpha -1}
\int^T_{s_{n+1}}f(u_{n+1})(u_{n+1})^{-\alpha}\\
&&\cdot \frac{ e(\delta_1 \vee \cdots \vee \delta_n \vee s_{n+1})
-e(\delta_1 \vee \cdots \vee \delta_{n}\vee u_{n+1})}
{(u_{n+1}-s_{n+1})^{1+\alpha}} du_{n+1}ds_{n+1}
\biggr] du_{i_j}\cdots du_{i_1}ds_n\cdots ds_1\\
&=&\frac{1}{\Gamma (\alpha)^{n+1}}\int^T_{t_1}\cdots \int^T_{t_{n+1}}
\frac{
\prod\limits^n_{i=1}\phi_f(s_i)s^{-\alpha}_i}
{\prod\limits^{n+1}_{i=1}(s_i- t_i)^{1-\alpha}}
\left[e(s_1\vee \cdots \vee s_{n+1})\phi_f(s_{n+1})s^{-\alpha}_{n+1}\right.\\
&& + \frac{\alpha}{\Gamma (1-\alpha)} \int^T_{s_{n+1}} \phi_f
(u_{n+1})u^{-\alpha}_{n+1} \frac{ e(s_1\vee \cdots \vee
s_{n+1})-e(s_1\vee \cdots \vee s_n \vee u_{n+1})}
{(u_{n+1}-s_{n+1})^{1+\alpha}}du_{n+1}\biggr]\\
&&\cdot ds_{n+1}\cdots ds_1\\
&& +\frac{1}{\Gamma(\alpha)^{1+n}} \sum^n_{j=1}
(\frac{\alpha}{\Gamma (1-\alpha)} )^j
\sum_{S_{j,n}}\int^T_{t_1}\cdots \int^T_{t_{n+1}} \frac{
(\prod\limits_{i\not\in \{i_1, \ldots , i_j, n+1\}} \phi_f (s_i)
(s^{-\alpha}_i)) \phi_f (s_{n+1})s^{-\alpha}_{n+1}}
{\prod^{n+1}_{i=1}(s_i-t_i)^{1+\alpha}}\\
&&\cdot \int^T_{s_{i_1}}\cdots \int^T_{s_{i_j}} \frac{
\prod\limits^{j}_{k=1}f(u_{i_k})u^{-\alpha}_{i_k}}
{\prod\limits^j_{k=1}(u_{i_k}-s_{i_k})^{1+\alpha}}
\sum^j_{k=0}(-1)^k \sum_{\delta \in \Delta_k^{
\widehat{s_{n+1}}}(u_{i_1},\ldots , u_{i_j})}
e(\delta_1 \vee \cdots \vee \delta_n \vee s_{n+1})\\
&&\cdot du_{i_j}\cdots du_{i_1}ds_{n+1}\cdots ds_1\\
&&+\frac{1}{\Gamma (\alpha)^{1+n}}\sum^{n+1}_{j=2}
(\frac{\alpha}{\Gamma
(1-\alpha)})^j\sum_{S_{j-1},n}\int^T_{t_1}\cdots \int^T_{t_{n+1}}
\frac{ \prod\limits_{i\not{\in} \{i_1, \ldots i_{j-1},
n+1\}}\phi_f (u_i) u^{-\alpha}_i }
{\prod\limits^{n+1}_{i=1}(s_i-t_i)^{1-\alpha}}\\
&&\cdot \int^T_{s_{i_1}}\cdots
\int^T_{s_{i_{j-1}}}\int^T_{s_{n+1}} \frac{
(\prod\limits^{j-1}_{k=1}f(u_{i_k}) u^{-\alpha}_{i_k})
f(u_{n+1})u^{-\alpha}_{n+1}}
{(\prod\limits^{j-1}_{k=1}(u_{i_k}-s_{i_k})^{1+\alpha})
(u_{n+1}-s_{n+1})^{1+\alpha}}
\sum^j_{k=0}(-1)^k\\
&&\cdot \sum_{\delta \in \Delta^s_k (u_{i_1},\ldots , u_{i_{j-1}},
u_{n+1})}e(\delta_1\vee \cdots \vee \delta_{n+1})du_{n+1}
du_{i_{j-1}}\cdots du_{u_{i_1}}ds_{n+1}\cdots ds_1,
\end{eqnarray*}
which proves that the result holds for $n+1$ whenever it is true for $n$. \hfill $\cajita$
\begin{lema}
\label{L4.4.1} Let $s\in [0,T]^n$, $n\geq 2$, $(i_1, \ldots, i_j)\in S_{j,n}$ and $(u_{i_1}, \ldots , u_{i_j})\in [s_{i_1},T]\times \cdots \times [s_{i_j}, T]$. Then
\begin{eqnarray*}
\lefteqn{
\sum^j_{k=0}(-1)^k \sum_{\delta \in \Delta^s_k (u_{i_1}, \ldots , u_{i_j})}
e(\delta_1\vee \cdots \vee \delta_n)}\\
&=&(e(s_1\vee \cdots \vee s_n)- e(u_{i_1}\wedge \cdots \wedge u_{i_j}))
1_{[s_1\vee \cdots \vee s_n, T]^j}(u_{i_1}, \ldots , u_{i_j}).
\end{eqnarray*}
\end{lema}
\noindent
{\it Proof:} The proof  follows from the fact that
\begin{eqnarray*}
&&\sum^{j+1}_{k=0}(-1)^k \sum_{\delta \in \Delta^s_k (u_{i_1},
\ldots, u_{i_j})} e(\delta_1 \vee \cdots \vee \delta_n)\\
&&=\sum^j_{k=0}(-1)^k \sum_{\delta \in
\Delta_k^{(\widehat{s_{i_{j_0}}})} (\widehat{u_{i_{j_0}}})}
[e(\delta_1\vee \cdots \vee \delta_{n-1}\vee s_{i_{j_0}})
-e(\delta_1 \vee \cdots \vee \delta_{n-1}\vee u_{i_{j_0}})],
\end{eqnarray*}
and induction on $j$. \hfill $\cajita$
\begin{lema}
\label{L4.4.2} Let $f\in \Lambda^\alpha_T$ be such that $\phi_f\in L^p([0,T])$ for some $p\in (2, 1/\alpha)$. Then
\begin{eqnarray}
\label{eq:4.4.1} \lefteqn{ ||\prod_{i\not\in \{i_1, \ldots,
i_j\}}\phi_f (s_i)
\int_{[s_1\vee \cdots \vee s_n,T]^j}}\nonumber\\
&&\cdot \frac{
(\prod\limits^j_{k=1}f(u_{i_k})u^{-\alpha}_{i_k}s^\alpha_{i_k})(e(s_1\vee \cdots \vee s_n)
-e(u_{i_1}\wedge \cdots \wedge u_{i_j}))}
{\prod\limits^j_{k=1}(u_{i_k}-s_{i_k})^{1+\alpha}}
du_{i_1}\cdots du_{i_j}||^2_{L^2([0,T]^n)}\nonumber\\
&\leq & jA^2 B^{2(j-1)}_p ||\phi_f||^{2(j-1)}_{L^p([0,T])}||\phi_f||^{2(n+1-j)}_{L^2([0,T])},
\end{eqnarray}
with
$$
A=\exp \left(\int^T_0
|a(s)|ds\right)||a||_{L^2([0,T])}\left(\int^1_0
\frac{d\theta}{(1-\theta)^{1-\alpha}\theta^{\frac{1}{2}+\alpha}}\right)
\frac{2T}{\Gamma(\alpha)}
$$
and
$$
B_p=\frac{T^{(p-2)/2p}}{\alpha}\left(\frac{p-2(1-\alpha
p)}{p-2}\right)^{\frac{p-2(1-\alpha p)}{2p}}C(p, p/(1-\alpha p)).
$$
\end{lema}

\noindent
{\it Proof:} Lemma \ref{L2.2.1} implies that the left--hand side of
(\ref{eq:4.4.1}) is bounded by
\begin{eqnarray*}
\lefteqn{
j!||\phi_f||^{2(n-j)}_{L^2([0,T])}
\exp ( 2\int^T_0 |a(s)|ds)||a||^2_{L^2([0,T])}}\\
&&\cdot \int^T_0\int^{s_j}_0\cdots \int^{s_2}_0\bigg((s_1\ldots
s_j)^\alpha \int_{[s_j,T]^j} \frac{
\prod\limits^j_{i=1}|f(u_i)|u_i^{-\alpha}}
{\prod\limits^j_{i=1}(u_i-s_i)^{1+\alpha}}
|u_j-s_j|^{1/2}du_1\ldots du_j\bigg)^2ds_1\ldots ds_j\\
&\leq & \frac{j}{\alpha^{2(j-1)}}\exp \left(2\int^T_0
|a(s)|ds\right)||a||^2_{L^2([0,T])}
\left(\int^1_0 \frac{(1-\theta)^{\alpha -1}}{\theta^{\frac{1}{2}+\alpha}}d\theta\right)^2 ||\phi_f||^{2(n-j)}_{L^2([0,T])}\\
&&\cdot \int^T_0 \left(\frac{1}{\Gamma (\alpha)}\int^T_{s_j}
\frac{|\phi_f(u)|}{(u-s_j)^{1/2}}du\right)^2\left[\int^{s_j}_0(s_j-s)^{-2\alpha}
\left(\frac{1}{\Gamma(\alpha)}\int^T_{s_j}\frac{|\phi_f(u)|}{(u-s)^{1-\alpha}}du\right)^2ds\right]^{j-1}ds_j.
\end{eqnarray*}
Hence, it is quite easy to finish the proof using Lemma \ref{L2.1.1}. \hfill $\cajita$
\begin{lema}
\label{L4.5.1} Let $f$ be as in Lemma \ref{L4.4.2}. Then
\begin{eqnarray*}
\lefteqn{
||f^{\otimes n}(t_1, \ldots, t_n)e (t_1\vee \ldots \vee t_n)
||_{(\Lambda^\alpha_T)^{\otimes n}}}\\
&\leq & C_H^{n/2}\biggl\{||\phi_f||^n_{L^2([0,T])}\exp (\int^T_0
|a(s)|ds)+ \\
&&\frac{nA\alpha}{\Gamma (1-\alpha)} \left(\frac{\alpha}{\Gamma
(1-\alpha)} B_p||\phi_f||_{L^p([0,T])}
+||\phi_f||_{L^2([0,T])}\right)^{n-1}\biggr\}.
\end{eqnarray*}
\end{lema}
\noindent
{\it Proof:} The proof is an immediate consequence of Lemmas \ref{L4.1.1},
\ref{L4.4.1} and \ref{L4.4.2}. \hfill $\cajita$
\begin{lema}
\label{L4.5.2} Let $f_n$ be given by (\ref{eq:3.2.1}).
Moreover assume that $\phi_b\in L^p([0,T])$, $p\in (2, 1/\alpha)$. Then
\begin{eqnarray*}
\lefteqn{\sum^\infty_{n=0}(n+1)!||\frac{1}{n+1}\sum^{n+1}_{k=1}1_{[0,t]}
(t_k)b(t_k)f^{t_k}_n(\hat{t}_k)||^2_{(\Lambda^\alpha_T)^{\otimes(n+1)}}}\\
&\leq & 2\exp\left(2 \int^t_0 |a(s)|ds
\right)\sum^\infty_{n=0}(n+1)!\left[\sum^n_{k=0}
\frac{1}{(n+1-k)!}||\eta_k||_{(\Lambda^\alpha_T)^{\otimes
k}}C_H^{(n+1-k)/2}
||\phi_{b1_{[0,t]}}||^{n+1-k}_{L^2([0,T])}\right]^2\\
&&+2\left(\frac{\alpha}{\Gamma (1-\alpha)}\right)^2
A^2\sum^\infty_{n=0}(n+1)!
\left[\sum^n_{k=0}\frac{1}{(n-k)!}C^{(n+1-k)/2}_H||\eta_k||_{(\Lambda^\alpha_T)^{\otimes k}}\right.\\
&&\quad \left.\cdot \left(\frac{\alpha}{\Gamma
(1-\alpha)}B_{\tilde{p}}||
\phi_{b1_{[0,t]}}||_{L^{\tilde{p}}([0,T])}+||\phi_{b1_{[0,t]}}||_{L^2([0,T])}\right)^{n-k}\right]^2,
\end{eqnarray*}
where $\tilde{p}\in (2,p)$ and $A, B_{\tilde{p}}$ are given in
Lemma \ref{L4.4.2}.
\end{lema}
\noindent
{\it Proof:} The definition of $f_n$ gives
\begin{eqnarray*}
\lefteqn{\frac{1}{n+1}\sum^{n+1}_{k=1}1_{[0,t]}(t_k)f^{t_k}_n(\hat{t}_k)b(t_k)}\\
&=&\sum^{n+1}_{k=1}\frac{1}{n+1}b(t_k)\exp (\int^{t_k}_0a(s)ds)\eta_n (\hat{t}_k)1_{[0,t]}(t_k)\\
&+&\sum^n_{j=1}\sum_{\Delta_{j+1, n+1}}\frac{(n-j)!}{(n+1)!}b^{\otimes(j+1)}
(t_{i_1},\ldots, t_{i_{j+1}})\eta_{n-j}(\hat{t}_{i_1},\ldots \hat{t}_{i_{j+1}})\\
&&\quad \cdot \exp(\int^{t_{i_{j+1}}}_0
a(s)ds)1_{\{t_{t_1}<\cdots <t_{j+1}<t\}}\\
&=&\sum^{n+1}_{k=1}\sum_{\Delta_{k, n+1}}
\frac{(n+1-k)!}{(n+1)!k!}(\prod^k_{j=1}b(t_{i_j})1_{[0,t]}
(t_{i_j}))\eta_{n+1-k}(\hat{t}_{i_1},\ldots, \hat{t}_{i_k})e
(t_{i_1}\vee \cdots \vee t_{i_k}).
\end{eqnarray*}
Hence Lemma \ref{L4.5.1} gives
\begin{eqnarray*}
\lefteqn{\sum^\infty_{n=0}(n+1)!||\frac{1}{n+1}\sum^{n+1}_{k=1}1_{[0,t]}
(t_k)b(t_k)f^{t_k}_n(\hat{t}_k)||^2_{(\Lambda^\alpha_T)^{\otimes (n+1)}}}\\
&\leq & \sum^\infty_{n=0}(n+1)!\left[\sum^{n+1}_{k=1}\frac{1}{k!}
||\eta_{n+1-k}||_{(\Lambda^\alpha_T)^{\otimes (n+1-k)}}
C^{k/2}_H\biggl\{||\phi_{b1_{[0,t]}}||^k_{L^2([0,T])}\exp (\int^t_0 |a(s)|ds)\biggr.\right.\\
&&\biggl.\left.\quad +\frac{A k \alpha}{\Gamma (1-\alpha)}
\left(\frac{\alpha}{\Gamma (1-\alpha)}B_{\tilde{p}} ||\phi_{b
1_{[0,t]}}||_{L^{\tilde{p}}([0,T])}+
||\phi_{b1_{[0,T]}}||_{L^2([0,T])}\right)^{k-1}\biggr\}\right]^2,
\end{eqnarray*}
and the result follows. \hfill $\cajita$
\begin{lema}
\label{L4.7.1} Assume that $\phi_b \in L^p ([0,T])$ for some $p\in (2, 1/\alpha)$. Then
\begin{eqnarray}
\label{eq:4.7.1}
\sum^\infty_{n=0}(n+1)!\left[\sum^n_{k=0}\frac{1}{(n+1-k)!}
||\eta_k||_{(\Lambda^\alpha_T)^{\otimes k}}C_H^{(n+1-k)/2}
||\phi_{b1_{[0,t]}}||^{n+1-k}_{L^2([0,T])}\right]^2\le \nonumber\\
\exp \left(1+C_H||\phi_{b1_{[0,t]}}||^2_{L^2([0,T])}\right)
\sum^\infty_{k=0}k!||\eta_k||^2_{(\Lambda^\alpha_T)^{\otimes k}}
(1+C_H||\phi_{b1_{[0,t]}}||^2_{L^2([0,T])})^k.
\end{eqnarray}
\end{lema}
\noindent
{\it Proof:} Observe that the left--hand side of (\ref{eq:4.7.1}) is bounded by
\begin{eqnarray*}
\lefteqn{
\sum^\infty_{n=0}(n+1)!\biggl[\sum^n_{k=0}\frac{k!}{(n+1-k)!}||\eta_k||^2_{(\Lambda^\alpha_T)^{\otimes k}}\biggr]
\biggl[\sum^n_{k=0}\frac{(C_H)^{n+1-k}}{k!(n+1-k)!}||\phi_{b1_{[0,t]}}||^{2(n+1-k)}_{L^2([0,T])}\biggr]}\\
&\leq & \sum^\infty_{k=0}k!||\eta_k||^2_{(\Lambda^\alpha_T)^{\otimes k}}(1+C_H||\phi_{b1_{[0,t])}}||^2_{L^2([0,T])} )^k
\sum^\infty_{n=k}\frac{(1+C_H||\phi_{b1_{[0,T]}}||^2_{L^2([0,T])})^{n+1-k}}{(n+1-k)!}.
\end{eqnarray*}
Thus, the proof is complete. \hfill $\cajita$
\begin{lema}
\label{L4.8.1}
Let $\tilde{p}\in (2,p)$. Then
\begin{eqnarray}
\label{eq:4.8.1}
&\sum\limits^\infty_{n=0}(n+1)!\biggl[\sum\limits^n_{k=0}\frac{C_H^{(n+1-k)/2}}{(n-k)!}
||\eta_k||_{(\Lambda^\alpha_T)^{\otimes
k}}\biggl(\frac{\alpha}{\Gamma (1-\alpha)}B_{\tilde{p}}||
\phi_{b1_{[0,t]}}||_{L^{\tilde{p}}([0,T])}+
||\phi_{b1_{[0,t]}}||_{L^2([0,T])}\biggr)^{n-k}\biggr]^2 \nonumber\\
&\leq C_H \exp (B_{H,\tilde{p},t})
\biggl\{\sum\limits^\infty_{k=0}(k+1)!||\eta_k||^2_{(\Lambda^\alpha_T)^{\otimes
k}}(B_{H,\tilde{p},t})^k +
\sum\limits^\infty_{k=0}k!||\eta_k||^2_{(\Lambda^\alpha_T)^{\otimes
k}} (B_{H,\tilde{p},t})^{k+1} \biggr\},
\end{eqnarray}
where
$$
B_{H,\tilde{p},t}=1+C_H \left(\frac{\alpha}{\Gamma (1-\alpha)}
B_{\tilde{p}}||\phi_{b1_{[0,t]}}
||_{L{^{\tilde{p}}([0,T])}}+||\phi_{b1_{[0,t]}}||_{L^2([0,T])}\right)^2.
$$
\end{lema}
\noindent
{\it Proof:}
Note that the left--hand side of (\ref{eq:4.8.1}) is dominated by
\begin{eqnarray*}
C_H \sum^\infty_{n=0}(n+1)\sum^n_{k=0}\frac{k!}{(n-k)!}||\eta_k
||^2_{(\Lambda^\alpha_T)^{\otimes k}}(B_{H,\tilde{p},t})^n=\\
C_H\sum^\infty_{k=0}k!||\eta_k||^2_{(\Lambda^\alpha_T)^{\otimes
k}}(B_{H,\tilde{p},t})^k \sum^\infty_{n=k}\frac{n+1}{(n-k)!}
(B_{H,\tilde{p},t})^{n-k}\leq \\
C_H \exp(B_{H,\tilde{p},t})
\sum^\infty_{k=0}k!||\eta_k||^2_{(\Lambda^\alpha_T)^{\otimes k}}
(B_{H,\tilde{p},t})^k \biggl(1+k +B_{H,\tilde{p},t}\biggr).
\end{eqnarray*}
Consequently, the result holds.  \hfill $\cajita$

\medskip

We will use the following results in the proof of Theorem \ref{T3.4.1}.

\begin{lema}
\label{L4.9.1} Let $p\in (2, 1/\alpha)$ and $\phi_b \in L^p
([0,T])$. Then for any $\delta <(\frac{p-2}{2p}\wedge \alpha
\wedge \frac{(1-\alpha p)}{p})$ there is a positive constant
$C_\delta$ such that
\begin{eqnarray*}
\lefteqn{
||(b1_{[0,t]})^{\otimes n}-(b1_{[0,s]})^{\otimes n}
||_{(\Lambda^\alpha_T)^{\otimes n}}}\\
&\leq & 2^{n-1}(C_\delta)^n ||\phi_b||^n_{L^p ([0,T])}(t-s)^\delta, \quad 0\leq s\leq t\leq T.
\end{eqnarray*}
\end{lema}
\noindent
{\it Proof:} We will use induction on $n$ to prove the result.

First assume that $n=1$. In this case, by (\ref{eq:2.3.1}), we have
\begin{eqnarray}
\label{eq:4.10.1}
\lefteqn{
||b1_{[0,t]}-b1_{[0,s]}||_{\Lambda^\alpha_T}}\nonumber\\
&\leq & C_H^{1/2}\biggl\{||\phi_{b}1_{]s,t]}||_{L^2([0,T])}+
\frac{\alpha}{\Gamma (1-\alpha)}\biggl(\int^s_0 \biggl(r^\alpha
\int^t_s \frac{b(u)u^{-\alpha}}{(u-r)^{1+\alpha}}du\biggr)^2dr\biggr)^{1/2}\nonumber\\
&&\quad + \frac{\alpha}{\Gamma (1-\alpha)}\biggl(\int^t_s
\biggl(r^\alpha \int^T_t \frac{b(u)u^{-\alpha}}{(u-r)^{1+\alpha}}du\biggr)^2dr\biggr)^{1/2}\biggr\}\nonumber\\
&=& C^{1/2}_H \biggl\{I_1 +\frac{\alpha}{\Gamma (1-\alpha)} I_2
+\frac{\alpha}{\Gamma (1-\alpha)}I_3\biggr\}.
\end{eqnarray}
It is clear that we have
\begin{equation}
\label{eq:4.10.2} I_1\leq (t-s)^{(p-2)/2p}||\phi_b
||_{L^p([0,T])}.
\end{equation}
Now, from Lemmas \ref{L2.1.1} and \ref{L2.2.1}, it follows
\begin{eqnarray}
\label{eq:4.10.3}
I_3 &\leq & \frac{1}{\Gamma (\alpha)\alpha}\biggl[\biggr(\int^t_s (t-r)^{-2\alpha}
\biggl(\int^T_r \frac{|\phi_b (u)|}{(u-r)^{1-\alpha}}du\biggr)^2dr\biggr]^{1/2}\nonumber\\
&\leq & \frac{C(p, p/(1-\alpha p))}{\alpha \Gamma(\alpha)}
\left(\frac{p-2(1-\alpha p)}{p-2}\right)^{(p-2(1-\alpha
p))/2p}||\phi_b ||_{L^p([0,T])}(t-s)^{(p-2)/2p}.
\end{eqnarray}

On the other hand, the Fubini's theorem and Lemma \ref{L2.2.1} lead to
\begin{eqnarray}
\label{eq:4.11.1}
\lefteqn{
1_{[0,s]}(r)r^\alpha \int^t_s\frac{|b(u)u^{-\alpha}|}{(u-r)^{1+\alpha}}du}\nonumber\\
&\leq & 1_{[0,s]}(r)\frac{r^\alpha}{\Gamma (\alpha)}\int^T_s |\phi_b (\theta)|\theta^{-\alpha}(\int^{t\wedge \theta}_s
\frac{du}{(u-r)^{1+\alpha}(\theta -u)^{1-\alpha}})d\theta\nonumber\\
&=&1_{[0,s]}(r)\frac{r^\alpha}{\Gamma (\alpha)}\biggl[\int^t_s |\phi_b (\theta)|\theta^{-\alpha}
(\int^\theta_s \frac{du}{(\theta -u)^{1-\alpha}(u-r)^{1+\alpha}})d\theta\nonumber\\
&&\quad + \int^T_t |\phi_b (\theta)|\theta^{-\alpha}(\int^t_s
\frac{du}{(u-r)^{1+\alpha}(\theta -u)^{1-\alpha}})d\theta \biggr]\nonumber\\
&\leq & 1_{[0,s]}(r)\biggl\{\frac{(t-s)^{\delta}}{\Gamma (\alpha)\alpha}(s-r)^{-\alpha}
\int^t_s \frac{|\phi_b (\theta)|}{(\theta-r)^{1-(\alpha -\delta)}}d\theta\nonumber\\
&&\quad +\frac{1}{\Gamma (\alpha)}(s-r)^{-\alpha-\delta}\int^t_s
\frac{1}{(u-r)^{1-\delta}}(\int^T_t \frac{|\phi_b (\theta)|}{(\theta-u)^{1-\alpha}}d\theta )du\biggr\}.
\end{eqnarray}
Therefore, combining (\ref{eq:4.10.1}) and (\ref{eq:4.11.1}), we get
\begin{eqnarray}
\label{eq:4.11.2}
I_2 &\leq & \frac{(t-s)^\delta}{\Gamma (\alpha)\alpha}(\int^s_0
(s-r)^{-2\alpha}(\int^T_r \frac{|\phi_b (\theta )|}{(\theta -r)^{1-(\alpha-\delta)}}d\theta )^2dr)^{1/2}\nonumber\\
&&\quad +\frac{1}{\Gamma (\alpha)} \biggl(\int^s_0
(s-r)^{-2(\alpha +\delta)}\biggl[\int^t_s
\frac{1}{(u-r)^{1-\delta}}(\int^T_u \frac{|\phi_b(\theta)|}{(\theta -u)^{1-\alpha}}d\theta)du\biggr]^2 dr
\biggr)^{1/2}\nonumber\\
&=& \frac{(t-s)^\delta}{\Gamma (\alpha)\alpha}I_{2,1}+
\frac{1}{\Gamma (\alpha)}I_{2,2}.
\end{eqnarray}
Observe that Lemma \ref{L2.1.1} gives
\begin{equation}
\label{eq:4.11.3}
I_{2,1}\leq [C(p, \frac{p}{1-p(\alpha -\delta)})(\frac{T^{1-2q\alpha}}{1-2q\alpha})^{1/q}]^{1/2}||\phi_b||_{L^p([0,T])},
\end{equation}
with $q=p(p-2(1-p(\alpha -\delta)))^{-1}$. Similarly, we have for $\delta <\frac{p-2}{p} \wedge \frac{2(1-\alpha p)}{p}$,
\begin{eqnarray*}
I_{2,2} &\leq & \biggl[\int^s_0 (s-r)^{-2(\alpha +\delta)}(\int^t_s \frac{du}{(u-r)^{1-\delta}})\\
&&\quad \cdot (\int^t_s \frac{1}{(u-r)^{1-\delta}}(\int^T_u
\frac{|\phi_b (\theta )|}{(\theta -u)^{1-\alpha}}d\theta)^2 du)dr \biggr]^{1/2}\\
&\leq & \frac{(t-s)^{\delta /2}}{\sqrt{\delta}}\biggl[\int^s_0 (s-r)^{-2(\alpha +\delta)}\int^T_r
\frac{1}{(u-r)^{1-\delta}}(\int^T_u \frac{|\phi_b (\theta)|}{(\theta -u)^{1-\alpha}}d\theta)^2dudr\biggr]^{1/2}\\
&\leq & \frac{(t-s)^{\delta /2}}{\sqrt{\delta}}
[\frac{T^{1-2(\alpha +\delta)\tilde{p}}}{1-2(\alpha +\delta)\tilde{p}}
]^{1/2\tilde{p}}C(p, \frac{p}{1-\alpha p})C(\frac{p}{2(1-\alpha p)},
\frac{p}{2(1-\alpha p)-\delta p} )||\phi_b||_{L^p ([0,T])},
\end{eqnarray*}
with $\tilde{p}=p(p-2(1-\alpha p)+\delta p)^{-1}$. Consequently,
from (\ref{eq:4.10.1}), (\ref{eq:4.10.2}),
(\ref{eq:4.10.3}),(\ref{eq:4.11.2}) and (\ref{eq:4.11.3}), we have
that the results holds for $n=1$.

Finally we can apply induction on $n$ to show the result is true for any $n$, due to Lemma \ref{L2.3.1} and the fact that
\begin{eqnarray*}
\lefteqn{||(b1_{[0,t]})^{\otimes n}-(b1_{[0,s]})^{\otimes n}||_{(\Lambda^\alpha_T)^{\otimes n}}}\\
&\leq & ||b1_{[0,s]}||_{\Lambda^\alpha_T}
||(b{1_{[0,t]}})^{\otimes (n-1)}-(b1_{[0,s]})^{\otimes (n-1)}||_{(\Lambda^\alpha_T)^{\otimes (n-1)}}\\
&& +||b1_{]s,t]}||_{\Lambda^\alpha_T}||(b1_{[0,t]})^{\otimes (n-1)}||_{(\Lambda^\alpha_T)^{\otimes (n-1)}}.
\end{eqnarray*}
\vglue -.5cm\hfill $\cajita$

An immediate consequence of Lemma \ref{L4.9.1} is the following.
\begin{corolario}
\label{C4.12.1}
Let $p, \delta$ and $\phi_b$ as in Lemma \ref{L4.9.1}. Then there is a constant $C$ such that
\begin{eqnarray*}
\lefteqn{\kern -3cm||\sum^n_{j=1}\sum_{\Delta_{j,n}}
\frac{(n-j)!}{n!j!}\eta_{n-j}(\hat{t}_{i_1},\ldots , \hat{t}_{i_j})
\biggl(\prod^j_{k=1}b(t_{i_k})1_{[0,t]}(t_{i_k})-\prod^j_{k=1}b(t_{i_k})1_{[0,s]}(t_{i_k})\biggr)
||_{(\Lambda^\alpha_T)^{\otimes n}}}\\
&\leq & (t-s)^\delta \sum^n_{j=1}\frac{1}{j!}||\eta_{n-j}||_{(\Lambda^\alpha_T)^{\otimes n-j}}
C^j||\phi_b||^j_{L^p ([0,T])}.
\end{eqnarray*}
\end{corolario}

\noindent\textsl{Acknowledgments}. The authors thank the Nucleus
Millennium \textsl{Information and Randomness} P01-005 for
support.

\end{document}